# Robust control of wheel slip using weighted fuzzy model

Mojtaba S. Zadeh, A. Akbari, M. Reza Zamani Behbahani

**Abstract:**

The control of wheel slip dynamics is one of the most critical areas of chassis control, since it is the basis for most critical of the main chassis control subsystems like braking control (ABS), traction control (TCS), and stability control (VSC-ESC). In this paper, a nonlinear adaptive controller is analytically designed for longitudinal wheel slip by using of the new method, based on Lyapunov theory mixed with weighted Fuzzy model. The effectiveness of the method is demonstrated by the experimental results and real time model simulations with different longitudinal speeds and different kinds of road surface conditions.

**Keywords:** NLLS, Optimum wheel slip, ESP, Recursive optimization, Antilock-braking system

# کنترل مقاوم لغزش وسیله نقلیه با اعمال الگوریتم فازی وزن دار

مجتبی شریف زاده[1] احمد اکبری[2] رضا زمانی[3]

[1] دانشجوی کارشناسی‌ارشدکنترل دانشگاه صنعتی سهند تبریز؛ دانشکده مهندسی برق و کامپیوتر دانشگاه تبریز؛   m_sharifzadeh@sut.ac.ir

[2] استادیار دانشکده مهندسی برق گروه کنترل دانشگاه صنعتی سهند تبریز؛   a.akbari@sut.ac.ir

[3] دانشجوی کارشناسی‌ارشدکنترل دانشگاه صنعتی سهند تبریز؛   mr_zamani@sut.ac.ir

چکیده

در این مقاله بەوسیله ترکیب روش‌های کنترل فازی و روش‌های کنترل تطبیقی و لیاپانوف یک قانون کنترل مقاوم با شرایط جاده‌ای متنوع طراحی شده‌است. سیستم‌های کنترل دینامیک خودرو، ایده نسبتا جدیدی در زمینه کنترل بلادرنگ خودرو می‌باشند که تحت عناوین مختلفی از قبیل ESC،VSC،AEB توسط شرکت‌های مختلف عرضه گشته‌اند. دراین میان استفاده از زیر سیستم ABS برای کنترل مستقل هرچرخ و نهایتا رسیدن به رفتار مطلوب خودرو اهمیت ویژه‌ای کسب کرده است.کارهای متفاوتی در زمینه کنترل لغزش به عنوان زیرمجموعه اساسی دراین بعد، صورت گرفته‌است که تحقیقات روی طراحی کنترل مقاوم نسبت به تغییرات دینامیکی وشرایط جاده‌ای ادامه دارد. نتایج شبیه‌سازی نشان می‌دهد که قانون کنترل پیشنهادی روی مقاومت بالایی نسبت به تغییرات شرایط دینامیکی وسیله و تخمین‌های صورت گرفته برای شرایط جاده‌ای مختلف دارا‌ست.

کلمات کلیدی: سیستم AEB، ESP، کنترل لغزش، ترمز ضدقفل

## ۱- مقدمه

اگرچه تعداد تصادفات در کشورهای دنیا درطول ده سال گذشته کاهش چندانی نداشته است، اما تعداد تلفات ناشی از این تصادفات بطور مداوم در حال کاهش است. در این بین، استفاده از سیستم های کنترلی جهت بهبود عملکرد دینامیکی خودرو و در نتیجه افزایش ایمنی آن به طور روزافزون در حال افزایش است و بسیاری از کارخانجات خودرو سعی در استفاده از سیستم هایی همچون ضدقفل، کنترل پایداری خودرو[1]، کنترل[2]، برنامه پایداری خودرو[2] ... دارند. لذا امروزه بحث کنترل و تخمین لغزش خودرو اهمیت ویژه ای پیدا کرده است. از آنجا که تشخیص کامل رفتار خودرو ، نیازمند نصب تعداد زیادی سنسور بر روی خودرو می‌باشد و از طرفی نصب این سنسورها قیمت خودرو را غیر اقتصادی می‌کند لذا از الگوریتم های تخمین و ریاضیات برای شناسایی رفتار خودرو استفاده می‌شود و کنترل سیستم بر این مبنا صورت می‌گیرد. نتایج حاصله نشان

می‌دهد[1] که نتایج تخمین تا حد قابل قبول به رفتار واقعی خودرو نزدیک است به گونه‌ای که کنترل حالت‌های ناشی از آن سیستم را پایدار می‌نماید شکل(۱).

اولین کاربرد سیستم ترمز ضدقفل (ABS[4]) در صنعت خودرو با الهام گیری از صنعت هواپیمایی توسط یک شرکت فرانسوی به اسم لینکلن در سال ۱۹۵۴صورت پذیرفت. هرچند کاربردهای نخست سیستم ضدقفل به جنگ جهانی دوم در جنگنده های B-47 باز میگردد[1] و برخی مقالات اشاره به ایده اولیه کاربرد آن در اوایل سال ۱۹۰۰ در سیستم قطار دارند. [2,3] اهمیت بالای ترمز ضدقفل و تاثیر آن با سیستم پردازنده آنالوگ تحقق نمی یافت. درهمین‌راستا شرکت های تویوتا و نیسان دستیابی به فناوری مبتنی برالکترونیک ضدقفل را اعلام کردند. سرانجام در سال ۱۹۷۸ شرکت بوش با همکاری شرکت دایملر-بنز موفق به ارائه سیستم ضدقفل گردید[2].

نسل جدید ترمزهای ضدقفل (ABS) و سیستم اضطراری ترمز اتوماتیک (AEB[4]) با تخمین هوشمند شرایط مختلف جاده ای و نهایتا رؤیت آتی

---

[3] Anti-Lock Braking System
[4] Automatic Emergency Braking

[1] Electronic Stability Control (Esc)
[2] Electronic Stability Prog. (Esp)

پارامترهای اصطکاکی اعم از نسبت لغزش[5] و ضریب اصطکاک[6] پارامترهای موردنیاز کنترل کننده را فراهم می آورد. به گونه ای که عمل ترمز در شرایط مختلف جاده ای اعم از برفی، آسفالت و .. با حداکثر بازدهی و نیروی ممکن صورت می پذیرد. گفتنی است در پی آن چرخ قفل نشده کنترل فرمان مطلوب تری خواهد داشت. شایان ذکر است، رؤیتگر اشاره شده در سایر قابلیت های نسل جدید خودروها اعم از کنترل پایداری الکترونیکی خودرو (ESC[7]) و سیستم کنترل لغزش خودرو (TCS[8]) نقش کاربردی ایفا می نماید [4،3].

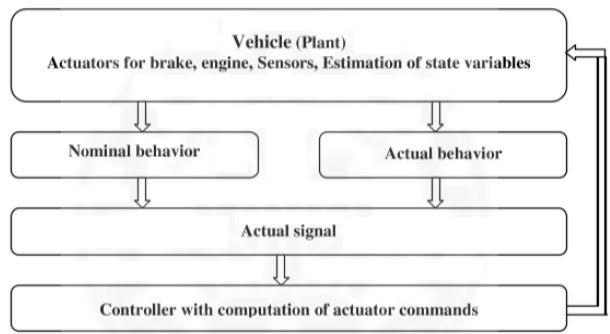

شکل۱. عملکرد کنترل پایداری خودرو

بخش بعدی به بیان روش‌های مدلسازی محوری مساله که مرجع اصلی تحقیقات امروزه را به خود اختصاص داده، پرداخته است. ابتدا مدلسازی خودرو در حالت Single-Corner بیان شده [5] و در ادامه مدلسازی اصطکاک بر حسب لغزش صورت گرفته است. [6]

دربخش۳. پایداری سیستم مورد بررسی قرار گرفته و در ادامه الگوریتم کنترلی بیان شده است . تابع تبدیل سیستم خطی شده در حالت کنترل لغزش محاسبه شده و نهایتا، پایدارسازی آن با کنترلر مورد بحث قرار گرفته است . در نهایت به بیان ایده کنترل تلفیقی فازی باتوجه به شرایط مختلف جاده‌ای پرداخته شده است .

درپایان بخش شبیه سازی ونتایج گزارشی از خروجی الگوریتم اعمال شده و مقایسه آن با نتایج گذشته می‌باشد که نشان می‌دهد الگوریتم بیان شده ضمن نداشتن حساسیت به نویز در مقابل تغییرات شرایط مختلف جاده‌ای مقاوم است.

## ۲- دینامیک و مدلسازی مساله

نخست به بیان دینامیک خودرو می پردازیم . جهت سهولت در بیان الگوریتم و راهکارهای کنترل و تخمین مدل ساده و استاندارد Single-Corner گزیده شده تا در صورت خروجی مطلوب مساله بتوان به مدل‌های ترکیبی نظیر

---

Double-Corner تعمیم داد [8،5].

$$\begin{cases} J\dot{w} = rF_x - T_b \\ m\dot{v} = -F_x \end{cases}$$

(۱)

در رابطه (۱) w[rad/s] سرعت زاویه ای چرخ است که توسط سنسور انکودر وسیله اندازه گیری می‌گردد. T_b [Nm] گشتاور ترمز ،v[rad/s] سرعت خطی وسیله نقلیه و r[m] ،m [Kg] و J[Kgm²] به ترتیب ممان اینرسی چرخ، شعاع چرخ و جرم Single-Corner می باشند.

$F_x$ نیروی اصطکاک چرخ و با رابطه (۲) قابل بیان می باشد:

$$F_x = F_z\mu(\lambda, \beta_t, v_r)$$

(۲)

که در آن

$F_z$ نیروی عمودی محل تقاطع تایر- جاده است.

$\beta_t، v_r$ به ترتیب ضرایب رابطه اصطکاک μ باتوجه به شرایط اصطکاک جاده ای و زاویه جانبی چرخ نسبت به صفحه عرضی می باشد [8،5]. $\lambda$ نسبت لغزش طولی بوده، $\lambda\epsilon[-1, 1]$ و به صورت رابطه

$$\lambda = (v - rw)/max\{rw, v\}$$

(۳)

بیان می گردد.

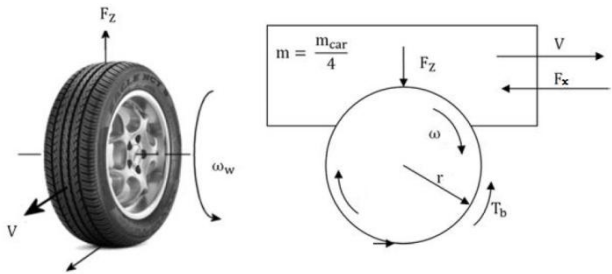

شکل ۲. مدل استاندارد Single-Corner

باتوجه به اینکه شرایط دینامیک مساله روی خط راست در نظر گرفته شده می توان ضریب زاویه $\beta$ را نادیده گرفت. دینامیک و شرایط مساله برای حالت ترمز (شتاب منفی) درنظر گرفته شده، لذا مخرج رابطه (۳) برابر $v$ درنظر گرفته شده است. نهایتا می توان رابطه دیفرانسیل (۴) را از آن استخراج نمود :

$$\lambda = -\frac{r}{v}\dot{w} + \frac{rw}{v^2}\dot{v}$$

(۴)

---

[5] Longitudinal Wheel Slip
[6] Friction coefficient

[7] Electronic Stability Control
[8] Traction Control System

نهایتا با جاگذاری روابط (۴) در (۳) و استفاده از رابطه (۲) دردینامیک سیستم، رابطه مهم (۵) نتیجه می گردد :

$$\begin{cases} \dot{\lambda} = -\frac{1}{v}\left(\frac{1-\lambda}{m} + \frac{r^2}{J}\right)F_z\mu(\lambda) + \frac{r}{vJ}T_b \\ m\dot{v} = -F_z\mu(\lambda) \end{cases}$$

(۵)

باتوجه به اینکه در سیستم موجود در رابطه(۵) تغییرات سرعت طولی سیستم (بیان شده با $v$) از تغییرات سرعت زاویه ای سیستم (بیان شده با $w$ یا $\lambda$) بسیار کمتر است، می توان از عبارت دوم رابطه (۵) صرف نظر کرد. [۸] شکل(۳) نمودار نیروی اصطکاک بر حسب لغزش را به ازای مقادیر متفاوت نیروی عمودی وارد بر چرخ نشان می‌دهد .

شکل ۴. نمودار ضریب اصطکاک برحسب لغزش- نتایج واقعی

برای عبارت ضریب اصطکاک بر حسب لغزش مدل های مختلفی بیان شده و قید شده اند،[۶] ازجمله مدل مجیک فرمولا و مدل بورکهارت که دو مدل معروف و اشاره شده در مقالات معتبر می باشند . [۱]،[۴] که در این مقاله از مدل معروف بورکهارت [۷] ¹⁰ بهره خواهیم جست.

$$\mu(\lambda; v_r) = V_{r1}\left(1 - e^{-\lambda v_{r2}}\right) - \lambda v_{r3}$$  (۶)

رابطه (۶) مدل بورکهارت را بازگو می‌سازد.

## ۳- تحلیل پایداری و طراحی کنترل کننده

حال یک سیستم غیرخطی بر حسب پارامتر $\lambda$ به صورت رابطه (۷) در اختیار داریم که $\mu(\lambda)$ خود تابعی برحسب پارامتر لغزش (رابطه ۶) می باشد.

$$\dot{\lambda} = -\frac{1}{v}\left(\frac{1-\lambda}{m} + \frac{r^2}{J}\right)F_z\mu(\lambda) + \frac{r}{vJ}T_b$$  (۷)

با درنظر گرفتن $v = \frac{rw}{1-\lambda}$ میتوان رابطه (۷) را به صورت زیر (۸) بازنویسی نمود :

$$\dot{\lambda} = \frac{\lambda - 1}{Jw}\left(\left(r + \frac{J(1-\lambda)}{rm}\right)F_z\mu(\lambda) - T_b\right)$$  (۸)

سیستم،رابطه (۸) یک سیستم غیرخطی به صورت (Autonomous) می باشد که ریشه های $\dot{\lambda} = 0$ نقاط تعادل سیستم را نتیجه می دهد. [۹] که با رسم سیستم فوق به ازای دو پارامتر $f(\lambda, \dot{\lambda})$ مشهود است.

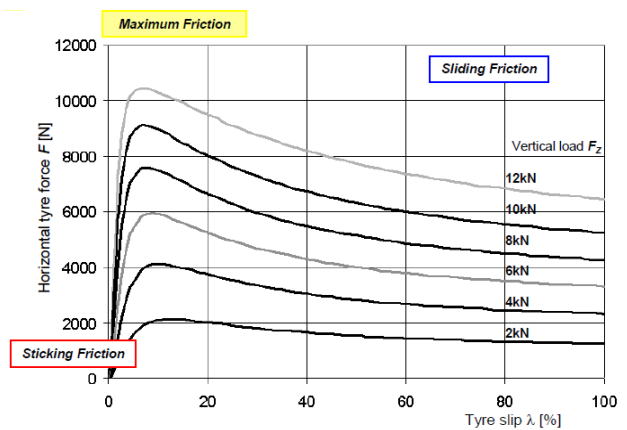

شکل ۳. نمودار نیروی اصطکاک بر حسب لغزش

بنابه رابطه (۲)، رابطه به صورت $\mu(\lambda; v_r) = F_x/F_z$ قابل بیان است . رابطه اصطکاک بر حسب ضریب لغزش در شرایط جاده- ای متفاوت به صورت نمودار شکل(۳) نمایان است . داده های مساله با شرایط آزمایشگاهی⁹ محاسبه گشته اند.[۸].

---

⁹ Experimental

می‌توان به صورت رابطه ۱۰ لحاظ کرد:

$$\begin{cases} \dot{z} = Az + \dfrac{r}{Jv}u, \\ y_{out} = z \end{cases}$$

(۱۰)

که در آن:

$$A = \frac{F_z}{v}\left[\frac{\mu(\bar{\lambda})}{m} - \dot{\mu}(\bar{\lambda})\left(\frac{1-\bar{\lambda}}{m} + \frac{r^2}{J}\right)\right]z,$$

جهت درنظر گرفتن سیستم (۸) به عنوان یک سیستم خطی سرعت طولی سیستم $v$ به عنوان یک سرعت ثابت فرض شده است. البته این مقدار طی بازه های گسسته توسط تخمین گر سیستم یا سیستم اندازه گیری سرعت طولی به روز رسانی می گردد که در ادامه بدان اشاره خواهیم داشت. با توجه به نکته (۱) و توضیحات فوق نهایتاً تابع تبدیل سیستم بر حسب خروجی لغزش و ورودی گشتاور ترمز به صورت رابطه (۱۱) به دست می آید . [۹].

$$G_\lambda(s) = \frac{\lambda}{T_b} = \frac{\dfrac{r}{Jv}}{s + \dfrac{F_z}{mv}\left[\dot{\mu}(\bar{\lambda})\left(1-\bar{\lambda}+\dfrac{mr^2}{J}\right) - \mu(\bar{\lambda})\right]}$$

(۱۱)

مشابها با اعمال رابطه تعریف شده برای شتاب منفی به صورت $\eta := -\dfrac{r\dot{w}}{g}$ داریم :

$$G_\eta(s) = \frac{\eta}{T_b} = \frac{\dfrac{r}{Jg}\left[s + \dfrac{F_z}{mv}\left(\dot{\mu}(\bar{\lambda})(1-\bar{\lambda}) - \mu(\bar{\lambda})\right)\right]}{s + \dfrac{F_z}{mv}\left[\dot{\mu}(\bar{\lambda})\left(1-\bar{\lambda}+\dfrac{mr^2}{J}\right) - \mu(\bar{\lambda})\right]}$$

(۱۲)

بنابراین سیستم خروجی شتاب منفی بر حسب ورودی گشتاور به صورت سیستم (۱۲) به دست می‌آید [۹].

**متدهای کنترل سیستم:**

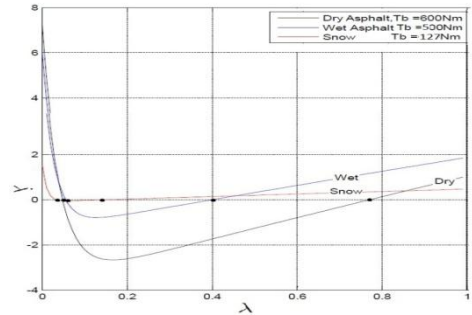

شکل ۵ ، نمودار $\dot{\lambda}$ بر حسب پارامتر $\lambda$

شکل (۵) نموداری از سیستم رابطه (۸) با دوپارامتر $f(\lambda, \dot{\lambda})$ به ازای شرایط مختلف جاده ای اعم از برفی، آسفالت خشک و آسفالت مرطوب می باشد. همان طور که از شکل (۵) واضح است این سیستم دارای دو نقطه تعادل است . مسیر تراژکتوری ها مشخص کننده این است که نقطه تعادل نخست پایدار مجانبی محلی و نقطه تعادل دوم ناپایدار می باشد [۹].

خطی سازی سیستم حول نقطه تعادل پایدار، سیستم خطی حول آن نقطه را نتیجه می دهد [۹].

**نکته ۱:**

فرض بر این سیستم غیرخطی به صورت (۷) داریم که نقاط تعادل آن $x_0 = (x_{10}, x_{20})$ می باشد.

$$\begin{cases} \dot{x}_1 = f_1(x_1, x_2) \\ \dot{x}_2 = f_2(x_1, x_2) \end{cases}$$

(۷)

سیستم خطی شده به صورت $\begin{cases} \dot{z}_1 = a_{11}z_1 + a_{12}z_2 \\ \dot{z}_2 = a_{21}z_1 + a_{22}z_2 \end{cases}$ نتیجه می‌شود که در آن داریم:

$$\dot{z} = Az, \quad A = \begin{bmatrix} a_{11} & a_{12} \\ a_{21} & a_{22} \end{bmatrix} = \begin{bmatrix} \dfrac{\partial f_1}{\partial x_1} & \dfrac{\partial f_1}{\partial x_2} \\ \dfrac{\partial f_2}{\partial x_1} & \dfrac{\partial f_2}{\partial x_2} \end{bmatrix}\Bigg|_{x = x_0}$$

(۸)

حال سیستم غیرخطی ۷ حالت خاصی از مساله فوق بوده و با خطی سازی حول نقطه تعادل پایدار $\bar{\lambda}$ و با درنظر گرفتن تغییر متغیرهای رابطه ۹:

$$x_1 = \lambda, x_2 = 0, y_{out} = z, T_b = u,$$

$$\dot{\mu}(\bar{\lambda}) = \frac{\partial \mu}{\partial \lambda}\Bigg|_{\lambda = \bar{\lambda}}$$

(۹)

در آن حدود ۱۰ میلی ثانیه و $w\_act = 70 rad/s$ درنظر گرفته می‌شود.

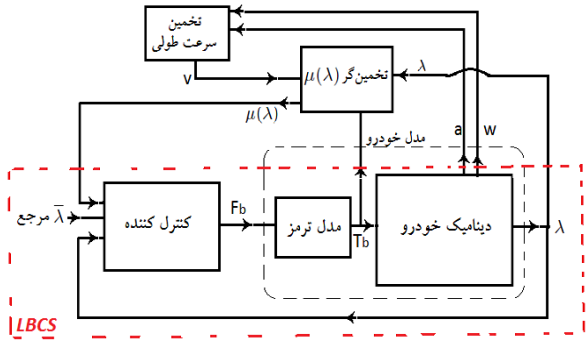

شکل ۷. ساختار کلی در کنترل پایداری لغزش

شکل (۷) مدل سیستم کنترل خطی ترمز را بازگو می‌سازد. اساس این سیستم بر مبنای کنترل لغزش خودرو و پایدارسازی آن در مقدار مرجع مورد نظر استوار است. با مراجعه به شکل (۱) بیشینه ضریب اصطکاک و نهایتا نیروی ترمزی بر حسب لغزش مورد تقاضا قابل بیان است و همین مقدار را می توان طی شناسایی به عنوان لغزش مرجع به کاربرد [۸].

روش تنظیم خودکار مبتنی بر فیدبک رله ای به فرم شکل (۸) می‌باشد. بروز رسانی پارامترهای PID براساس روابط بهره و پریود بحرانی (۱۷) محقق می‌گردد [۱۰].

$$\arg G(jw_u) = -\pi \ , \qquad |G(jw_u)| = \frac{1}{k_u} \ , \quad w_u = \frac{2\pi}{p_u}$$
$$(۱۷)$$

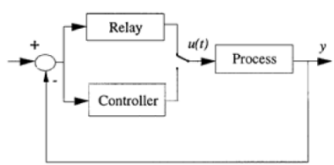

شکل ۸. کنترل PID خود تنظیم مبتنی بر فیدبک رله‌ای

در ذیل به بیان کنترل سیستم به روش لیاپانف پرداخته شده‌است. بنابر رابطه(۸) می‌توان نوشت .

$$\dot{\lambda} = -\frac{1-\lambda}{jw}\left(\left(\Psi(\lambda)\right) - T_b\right) \qquad (۱۸)$$
$$\Psi(\lambda) = \left(r + \frac{J(1-\lambda)}{rm}\right)F_z\mu(\lambda)$$

شکل (۶) نمودار اصطکاک بر حسب لغزش را نشان می‌دهد و به بیان پتانسیل اصطکاک پرداخته است. همان طور که از شکل (۷) بر می آید لغزش در یک مقدار بهینه دارای پیک اصطکاک می‌باشد. ایده کنترلی بر این پایه استوار است که با کنترل لغزش در مقدار بهینه اصطکاک سیستم را در حداکثر نگهداریم.

باتوجه به رابطه تابع تبدیل لغزش سیستم (۱۱)، معادله مشخصه سیستم به صورت (۱۵) قابل بیان است.

$$\varphi_\eta(s) = s + \frac{1}{v}\left[\frac{\dot{\mu}(\bar{\lambda})F_z}{m}\left(1 - \lambda + \frac{mr^2}{J}\right) + K\frac{r}{J}\right] - \frac{\mu(\bar{\lambda})F_z}{mv}$$
$$(۱۵)$$

که نهایتا شرط پایداری سیستم با بررسی وضعیت ریشه های رابطه (۱۵) به صورت رابطه(۱۶) بدست می آید.

$$K > -\frac{\dot{\mu}(\bar{\lambda})F_zJ}{mr}\left((1 - \bar{\lambda}) + \frac{mr^2}{J}\right) + \frac{\mu(\bar{\lambda})F_zJ}{mr}$$
$$(۱۶)$$

از رابطه (۱۶) استنباط می گردد کنترل پایداری لغزش سیستم نیازمند تخمین ضریب اصطکاک بر حسب لغزش خواهد بود. لذا همان‌گونه که در سیستم شکل (۸) مشخص شده محاسبه یا برآورد ضرایب مجهول فرم بور کهارت (۶) بحثی جداگانه می طلبد [۸].

مزیت کنترل به روش پایداری لغزش، امکان کنترل مقاوم نسبت به شــرایط مختلف جاده ای و نیز در سرعت های مختلف می باشـد [۱۱]. هرچند ایراد کار عامل نویز در سرعت های پایین می باشد [۱۱].

دینامیک حلقه بسته کنترلر سرو کنترل ترمز الکترومکانیکی درنظر گرفته شده EMB دارای تابع تبدیل سیستم مرتبه اول دارای تاخیر،باورودی نیرو در خروجی گشتاور ترمز به فرم $G_b = \frac{w_{act}}{s+w_{act}}e^{-sT}$ می باشد که زمان تاخیر

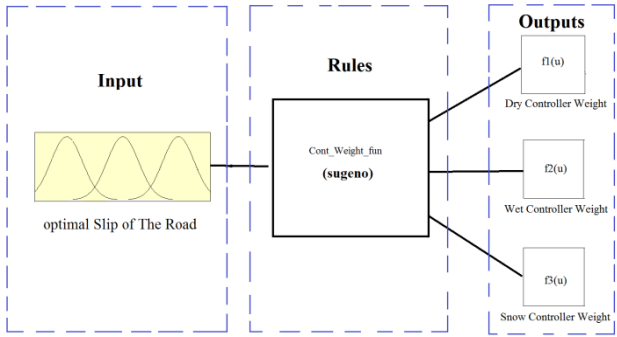

شکل ۱۰. ساختار فازی ارائه شده در مساله

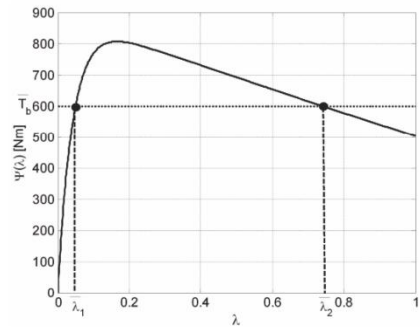

شکل ۹. نمودار گشتاور ترمز برحسب لغزش

در شکل (۹) با توجه به مسیر تراژکتوری ها مشخص است که نقطه تعادل نخست پایدار میجانبی محلی و نقطه تعادل دوم ناپایدار می باشد. تابع لیاپانف به فرم زیر پایداری سیستم را تامین میکند.

$$W(\lambda,\theta) = -\lambda + (\overline{\lambda} - 1)\ln(1-\lambda) + \varepsilon(\theta) + c \quad (۱۹)$$

که در آن

$$\varepsilon(\theta) := \ln(\text{T}_{\text{bMax}} - \theta)^{\frac{\tau}{k_\lambda}} - \ln(\theta - \text{T}_{\text{bMin}})^{\frac{\tau+1}{k_\lambda}}$$

و $\tau = \dfrac{\overline{\theta} - \text{T}_{\text{bMax}}}{\text{T}_{\text{bMax}} - \text{T}_{\text{bMin}}} < 0 \,,\ \tau < -1$ میباشد.

شکل(۱۰) ساختار فازی[۱۱] مساله را بیان مینماید. ورودی لغزش بهینه[۱۲] به ازای شرایط مختلف جادهای مقداری بازهای است که میتواند به عنوان ورودی مساله بکار گرفته میشود. در شکل (۱۱) عملکرد فازی به صورت گروهها بیان شده است. نمودار نمونه نشانداده مثالی است برای شرایط جادهای مابین مرطوب و برفی. به طوری که موجب میگردد کنترلرهای مربوط به آنها با وزن محاسبه شده از الگوریتم فازی عمل نمایند.

$\theta$ درنظر گرفته شده در رابطه فوق میباشد و $T_b$ مساله و نهایتا قانون کنترل به صورت رابطه (۲۰) بیان میگردد.

$$w\dot{\theta} = k_\lambda \frac{1}{J} \left(\lambda - \overline{\lambda}\right)(\theta - \text{T}_{\text{bMax}})(\theta - \text{T}_{\text{bMin}}) \quad (۲۰)$$

$\text{J}[\text{Kgm}^2]$ ممان اینرسی چرخ میباشد. تفاوت ساختاری ارائهشده با شکل ۸ در این است که کنترلر و دینامیک ترمز یکجا درنظر گرفته شده و خروجی کنترلر به صورت گشتاور ترمزی (Tb) ارائه میگردد و همانگونه که در رابطه فوق اشاره گشت، بهروز رسانی میشود.

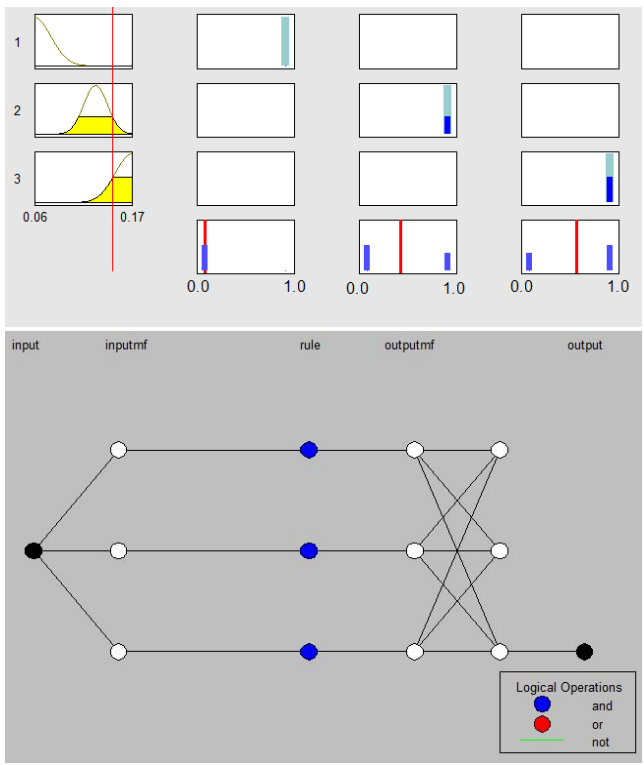

شکل ۱۱ فلوچارت و نمودار نمونه الگوریتم فازی وزن دار

هردو روش کنترل ارائه شده پایداری سیستم را در با تغییرات سرعت محقق میسازند. هرچند خروجی روش لیاپانف نسبت به تغییرات نویز مطلوبتر است. مساله اصلی مقاوم بودن کنترل کننده نسبت به تغییرات شرایط جادهای است. جهت این مهم الگوریتم جدیدی را جهت تلفیق با یکی از روشهای کنترلی فوق به صورت زیر تعریف میکنیم. الگوریتم فازی توصیف شده در موازات قانون کنترلی فوق میباشد به گونهای که برای شرایط جادهای وزنهایی را در نظر گرفته میشود که به عنوان خروجی الگوریتم فازی مورد بحث قرار میگیرد.

---

[12] Optimal Slip or $\lambda_{opt}$

[11] Fuzzy Structure

شکل(۱۱) ضرایب کنترلی اعمالی در کنترل‌کننده‌های پارالل را نشان می‌دهد که از الگوریتم فازی استنتاج شده‌اند. همان‌طور که در شکل(۴) بیان شد، شرایط جاده‌ای صرف به سه شرط اصلی فوق بسته نبوده و شرایط ترکیبی اعم از سنگفرش[13] وغیره وجود دارد. با توجه به اینکه نتایج شبیه‌سازی، تطبیق‌پذیری شرایط مختلف تست‌شده را بازگو می‌نماید و صرفا جهت گویا بودن الگوریتم بیان‌شده و بدون اینکه در کلیت مساله اشکالی وارد گردد سه حالت مختلف اصلی جاده‌ای درمساله آورده شده‌اند.

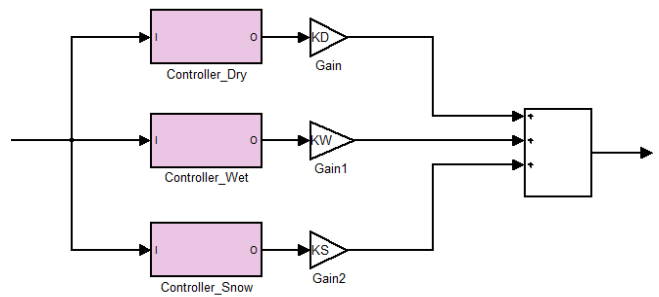

شکل۱۳. نمودار لغزش کنترل شده بر حسب زمان و شرایط متفاوت دینامیکی و جاده‌ای

شکل۱۲. ساختار کنترل وزن داده شده

## ۴– شبیه‌سازی و نتایج :

شبیه‌سازی سیستم فوق به کمک محیط سیمولینک و دستورات نرم افزار متلب و جعبه‌ابزار فازی انجام شده است. نمودار لغزش کنترل شده بر حسب زمان در شرایط متفاوت دینامیکی و جاده‌ای در نمودار شکل(۱۳) بیان شده است. همان طور که مشخص است نمونه شبیه‌سازی شده با اعمال ورودی لغزش مرجع از ثانیه یکم و با درنظرگیری اغتشاش واحد بعد از ثانیه دوم در شرایط مختلف جاده‌ای در سرعت‌های مختلف انتخابی اجرا گشته است.

خروجی شبیه نشان می‌دهد با وجود تغییرات جاده‌ای و در سرعت‌های متفاوت مقدار لغزش به مقدار مطلوب همگرا شده است . نتایج حاکی از آن است حساسیت کنترل کننده در سرعت‌های پایین تر از v=5m/s به نویز افزایش می‌یابد و القای ناپایداری و پایدار نوسانی در سیستم را دربی دارد.هرچند از دید عملی این نگرانی در سرعت‌های متمایل به صفر مطرح نیست. در سرعت‌های بیش‌تر از مقدار مذکور خروجی مطلوب و پایدار رویت می‌گردد. همچنین در نتایج شبیه سازی ملاحظه می‌گردد مقاوم روش مذکور نسبت به شرایط متنوع جاده ای بالاست. امروزه جهت کنترل ترمز روش های غیرخطی و مقاوم متنوعی مطرح شده است[۱،۳،۴،۶]. از مزایای روش های کنترلی غیرخطی به روش‌های مبتنی بر PID ذکرشده دراین است که حساسیت به نویز کمتری در خروجی احساس می‌گردد [۱۱]. حال نقطه قوت الگوریتم تلفیق داده‌شده دراین است که نشان می‌دهد با اعمال قوانین فازی و وزن دادن به کنترلر مساله می‌توان کنترل را نسبت به شرایط متنوع جاده‌ای مقاوم کرد.

---

[13] Cobblestone

مراجع